\documentclass{amsart}

\usepackage{amsmath,amssymb,amsthm}
\usepackage[all]{xy}

\DeclareMathOperator{\Diff}{Diff}
\DeclareMathOperator{\Aut}{Aut}
\DeclareMathOperator{\Lie}{Lie}
\DeclareMathOperator{\Hom}{Hom}
\DeclareMathOperator{\SL}{SL}
\DeclareMathOperator{\GL}{GL}
\DeclareMathOperator{\Or}{O}

\DeclareMathOperator{\hol}{hol}
\DeclareMathOperator{\id}{id}

\DeclareMathOperator{\pr}{pr}

\DeclareMathOperator{\Ad}{Ad}

\newtheorem{Theorem}{Theorem}
\newtheorem{Lemma}[Theorem]{Lemma}
\newtheorem{Proposition}[Theorem]{Proposition}
\newtheorem{Corollary}[Theorem]{Corollary}

\begin{document}
\title{Minimizability of developable Riemannian foliations}
\author{Hiraku Nozawa}
\address{Unit\'{e} de Math\'{e}matiques Pures et Appliqu\'{e}es, \'{E}cole Normale Sup\'{e}rieure de Lyon, \newline
46, all\'{e}e d'Italie 69364 Lyon Cedex 07, France.}
\email{nozawahiraku@06.alumni.u-tokyo.ac.jp}

\begin{abstract}
Let $(M,\mathcal{F})$ be a closed manifold with a Riemannian foliation. We show that the secondary characteristic classes of the Molino's commuting sheaf of $(M,\mathcal{F})$ vanish if $(M,\mathcal{F})$ is developable and $\pi_{1}M$ is of polynomial growth. By theorems of \'{A}lvarez L\'{o}pez in \cite{Alvarez Lopez} and \cite{Alvarez Lopez 2}, our result implies that $(M,\mathcal{F})$ is minimizable under the same conditions. As a corollary, we show that $(M,\mathcal{F})$ is minimizable if $\mathcal{F}$ is of codimension $2$ and $\pi_{1}M$ is of polynomial growth.
\end{abstract}

\keywords{Riemannian foliation, Taut foliations, Secondary characteristic classes}

\maketitle

\section{Introduction}

A foliated manifold $(M,\mathcal{F})$ is minimizable if there exists a Riemannian metric $g$ on $M$ such that every leaf of $\mathcal{F}$ is a minimal submanifold of $(M,g)$. The minimizability of general foliations is characterized by dynamical tools; for example, foliation cycles \cite{Sullivan} or holonomy pseudogroups \cite{Haefliger}. On the other hand, the minimizability of Riemannian foliations has a strong relation with the topology of manifolds. Ghys \cite{Ghys} showed that every Riemannian foliation on a closed simply connected manifold is minimizable. \'{A}lvarez L\'{o}pez \cite{Alvarez Lopez} defined a cohomology class of degree $1$ called the \'{A}lvarez class whose triviality is equivalent to the minimizability of $(M,\mathcal{F})$. As a corollary, he showed that every Riemannian foliation on a closed manifold whose first Betti number is zero is minimizable. Masa \cite{Masa} characterized the minimizability of Riemannian foliations in terms of the basic cohomology. In this article, we show that the secondary characteristic classes of the Molino's commuting sheaf of a closed manifold $M$ with a Riemannian foliation $\mathcal{F}$ vanish under some topological conditions of $(M,\mathcal{F})$. The Molino's commuting sheaf of $(M,\mathcal{F})$ is a locally constant sheaf defined by Molino \cite{Molino} for a closed manifold with a Riemannian foliation. According to \'{A}lvarez L\'{o}pez \cite{Alvarez Lopez 2}, the \'{A}lvarez class of $(M,\mathcal{F})$ is equal to the secondary characteristic classes of degree $1$ of the Molino's commuting sheaf of $(M,\mathcal{F})$ up to multiplication of a non-zero constant. As a corollary, we show that a developable Riemannian foliation on a closed manifold with fundamental group of polynomial growth is minimizable. As an application of this result, we show that a Riemannian foliation of codimension $2$ on a closed manifold  with fundamental group of polynomial growth is minimizable.

A foliated manifold $(M,\mathcal{F})$ is developable if there exists a covering space $\pi \colon \tilde{M} \longrightarrow M$ such that the leaves of $\pi^{*}\mathcal{F}$ are fibers of a submersion. Our main result in this article is the following: Let $(M,\mathcal{F})$ be a closed manifold with a Riemannian foliation.
\begin{Theorem}\label{Theorem : Vanishing}
The secondary characteristic classes of the Molino's commuting sheaf of $(M,\mathcal{F})$ vanish if $\mathcal{F}$ is developable and $\pi_{1}M$ is of polynomial growth.
\end{Theorem}
By Theorem \ref{Theorem : Vanishing} and theorems of \'{A}lvarez L\'{o}pez in \cite{Alvarez Lopez} and \cite{Alvarez Lopez 2}, we have
\begin{Corollary}\label{Corollary : Minimizability}
A developable Riemannian foliation on a closed manifold with fundamental group of polynomial growth is minimizable.
\end{Corollary}
\noindent Each of the assumption of the developability and the polynomial growth of fundamental groups in Corollary \ref{Corollary : Minimizability} is essential. The Carri\`{e}re's examples are non-minimizable developable Riemannian flows on $3$-dimensional manifolds with the fundamental group of exponential growth (see \cite{Carriere} and Section \ref{Subsection : Carriere}). Example \ref{Example} provides a non-minimizable non-developable Riemannian flow on a closed manifold with abelian fundamental group. Example \ref{Example} is the first example of a non-minimizable Riemannian foliation on a closed manifold with fundamental group of polynomial growth in the literature as far as the author knows.

As we will see in Proposition \ref{Lemma : Codimension 3}, every transversely parallelizable foliation of codimension less than $4$ is compact or developable. Since a closed manifold with a compact Riemannian foliation is minimizable (see Proposition \ref{Proposition : Compact}), we have
\begin{Theorem}\label{Theorem : Codimension 3}
A transversely parallelizable foliation of codimension $q \leq 3$ on a closed manifold  with fundamental group of polynomial growth is minimizable.
\end{Theorem}

By Theorem \ref{Theorem : Codimension 3}, we have
\begin{Corollary}\label{Corollary : Codimension 2}
A Riemannian foliation of codimension $2$ on a closed manifold with fundamental group of polynomial growth is minimizable.
\end{Corollary}
\noindent Since the Carri\`{e}re's examples are of codimension $2$, the assumption of polynomial growth of fundamental groups is essential. Note that \'{A}lvarez L\'{o}pez \cite{Alvarez Lopez} showed that every Riemannian foliation of codimension $1$ on a closed manifold is minimizable.

The results in this article were obtained during the stay of the author at the University of Santiago de Compostela in the spring of 2009. The author would like to express his gratitude to Jes\'{u}s Antonio \'{A}lvarez L\'{o}pez for his invitation, great hospitality and inspiring discussion with him. 

\section{Unipotency of the holonomy homomorphism of the Molino's commuting sheaves}

Let $(M,\mathcal{F})$ be a closed manifold with a transversely parallelizable foliation. By the Molino's structure theorem \cite{Molino}, $M$ is the total space of a smooth fiber bundle $p \colon M \longrightarrow W$ such that each fiber of $p$ is the closure of a leaf of $\mathcal{F}$. This $p$ is called the basic fibration of $(M,\mathcal{F})$. Let $F$ a fiber of $p$. Then $p$ is a smooth $(F,\mathcal{F}|_{F})$-bundle by the Molino's structure theorem \cite{Molino}. Let $x_{0}$ be a point in $F$. Let
\begin{equation}
\partial \colon \pi_{2}(W,p(x_{0})) \longrightarrow \pi_{1}(F,x_{0})
\end{equation}
be the connecting homomorphism in the homotopy exact sequence of the fiber bundle $p$. $(F,\mathcal{F}|_{F})$ is a $G$-Lie foliation with dense leaves for a connected simply connected Lie group $G$ by Theorem 4.2 in \cite{Molino}. Let
\begin{equation}
\hol_{F} \colon \pi_{1}(F,x_{0}) \longrightarrow G
\end{equation}
be the holonomy homomorphism of the Lie foliation $(F,\mathcal{F}|_{F})$. Let $H$ be the closure of the image of $\hol_{F} \circ \partial$ in $G$. Note that since $\pi_{2}(W,p(x_{0}))$ is abelian, $H$ is abelian. We define $\Aut(G,H)$ by
\begin{equation}
\Aut(G,H) = \{ a \in \Aut(G) \ | \ a(H) = H \}.
\end{equation}

Let $\Diff(F,x_{0},\mathcal{F}|_{F})$ be the group of diffeomorphisms of $F$ which fix $x_{0}$ and preserve $\mathcal{F}|_{F}$. We recall the definition of the holonomy homomorphism $h_{p} \colon \pi_{1}(M,x_{0}) \longrightarrow \pi_{0}(\Diff(F,x_{0},\mathcal{F}|_{F}))$ of the $(F,\mathcal{F}|_{F})$-bundle with base points. Let $\gamma \colon S^1 \longrightarrow M$ be a smooth map such that $\gamma(1)=x_{0}$. We consider the bundle $(\pi \circ \gamma)^{*}\pi \colon (\pi \circ \gamma)^{*}M \longrightarrow S^1$ obtained by pulling back $\pi$ by $\pi \circ \gamma$. Then $\gamma$ gives a fiberwise base point $b \colon S^1 \longrightarrow (\pi \circ \gamma)^{*}M$ of $(\pi \circ \gamma)^{*}\pi$ defined by $b(t)=(t,\gamma(t))$. Note that $(\pi \circ \gamma)^{*}M=\{(t,x) \in S^1 \times M \ | \ \pi \circ \gamma(t) = \pi(x)\}$. Let $\mathcal{G}$ be the foliation on $(\pi \circ \gamma)^{*}M$ obtained by pulling back $\mathcal{F}$ by the map $(\pi \circ \gamma)^{*}M \longrightarrow M$ induced by $\pi \circ \gamma$. Let $\mathcal{G}_{0}$ be a foliation on $([0,1] \times F)/(0,f(x)) \sim (1,x)$ induced from the foliation with leaves $\{t\} \times L$, for $t \in [0,1]$ and $L \in \mathcal{F}$, on $[0,1] \times F$. There exists an element $f_{\gamma}$ of $\Diff(F,x_{0},\mathcal{F}|_{F})$ and an isomorphism $\chi \colon (\pi \circ \gamma)^{*}M \longrightarrow ([0,1] \times F)/(0,f_{\gamma}(x)) \sim (1,x)$ as an $F$-bundle over $S^1$ such that $\chi^{*}\mathcal{G}_{0}=\mathcal{G}$ and $\chi$ maps the fiberwise base point $\{(t,\gamma(t)) \in (\pi \circ \gamma)^{*}M \ | \ t \in S^1 \}$ to $\{ (t,x_{0}) \in [0,1] \times F \ | \ t \in [0,1] \}/(0,f_{\gamma}(x)) \sim (1,x)$. Note that we regard $([0,1] \times F)/(0,f_{\gamma}(x)) \sim (1,x)$ as an $F$-bundle over $S^1$ by the map $\pi_{1}$ defined by $\pi_{1}([t,x])=[t]$. The isotopy class $[f_{\gamma}]$ of $f_{\gamma}$ in $\Diff(F,x_{0},\mathcal{F}|_{F})$ is well-defined for the element $[\gamma]$ of $\pi_{1}(M,x_{0})$ determined by $\gamma$. Thus we have a homomorphism
\begin{equation}
\begin{array}{cccc}
h_{p} \colon & \pi_{1}(M,x_{0}) & \longrightarrow & \pi_{0}(\Diff(F,x_{0},\mathcal{F}|_{F})) \\
             & [\gamma]           & \longmapsto     & [f_{\gamma}].
\end{array}
\end{equation}
We define 
\begin{equation}
\pi_{1} \colon \pi_{0}(\Diff(F,x_{0},\mathcal{F}|_{F})) \longrightarrow \Aut(\pi_{1}(F,x_{0}))
\end{equation}
by $\pi_{1}([f])([\delta]) = [f \circ \delta]$ for $[f]$ in $\pi_{0}(\Diff(F,x_{0},\mathcal{F}|_{F}))$ and $[\delta]$ in $\pi_{1}(F,x_{0})$. We define
\begin{equation}
\rho_{0} \colon \pi_{1}(M,x_{0}) \longrightarrow \Aut(\pi_{1}(F,x_{0}))
\end{equation}
by $\rho_{0} = \pi_{1} \circ h_{p}$. 

\begin{Lemma}\label{Lemma : Pi 2}
The subgroup $\partial (\pi_{2}(W,p(x_{0})))$ of $\pi_{1}(F,x_{0})$ is preserved by the action $\rho_{0}$.
\end{Lemma}

\begin{proof}
Let $\gamma$ be a closed path in $F$. Let $\phi \colon D^2 \longrightarrow M$ be a $2$-disk such that $\phi|_{\partial D^2}=\gamma$. For $\delta$ in $\pi_{1}(M,x_{0})$, there exists an annulus $\psi \colon [0,1] \times S^1 \longrightarrow M$ such that $\psi|_{\{0\} \times S^1}=\gamma$ and $\psi|_{\{1\} \times S^1} = \rho_{0}(\delta)(\gamma)$. Hence $\rho_{0}(\delta)(\gamma)$ is bounded by the $2$-disk obtained by gluing $\phi$ and $\psi$ along $\gamma$. Thus $\pi_{1}([f])(\gamma)$ is an element of $\partial (\pi_{2}(W,p(x_{0})))$. 
\end{proof}

We briefly recall the definition of the Molino's commuting sheaf of transversely parallelizable foliation $(M,\mathcal{F})$. For an open set $U$ in $M$, $C_{b}^{\infty}(TU/T(\mathcal{F}|_{U}))$ is the space of transverse fields of $(U,\mathcal{F}|_{U})$ defined by
\begin{equation}
C_{b}^{\infty}(TU/T(\mathcal{F}|_{U})) = \{ Y \in C^{\infty}(TU/T(\mathcal{F}|_{U})) \ | \ [Z,Y]=0, \forall Z \in C^{\infty}(T(\mathcal{F}|_{U})) \}.
\end{equation}
We define a vector space $\mathcal{C}(U)$ as the space of transverse fields which commute with global transverse fields as follows: 
\begin{equation}
\mathcal{C}(U) = \{ X \in C^{\infty}_{b}(TU/T(\mathcal{F}|_{U})) \ | \ [X,Y|_{U}]=0, \forall Y \in C_{b}^{\infty}(TM/T\mathcal{F}) \}.
\end{equation}
Note that Lie brackets on vector fields on $M$ induces Lie brackets $[\cdot,\cdot] \colon C^{\infty}(T\mathcal{F}) \times C^{\infty}(TM/T\mathcal{F}) \longrightarrow C^{\infty}(TM/T\mathcal{F})$ and $[\cdot,\cdot] \colon C_{b}^{\infty}(TM/T\mathcal{F}) \times C_{b}^{\infty}(TM/T\mathcal{F}) \longrightarrow C_{b}^{\infty}(TM/T\mathcal{F})$. Since $\mathcal{C}(U)$ is closed under the Lie bracket, $\mathcal{C}$ is a presheaf of Lie algebras. The sheafication of $\mathcal{C}$ is called the Molino's commuting sheaf of $(M,\mathcal{F})$.

We refer to pages 125--130 of Molino \cite{Molino} for the basic properties of Molino's commuting sheaf of transversely parallelizable foliations. Let $\mathcal{C}$ be the Molino's commuting sheaf of $(M,\mathcal{F})$. By Theorem 4.3 of Molino \cite{Molino}, $\mathcal{C}$ is locally constant and the stalk $\mathcal{C}_{x_{0}}$ of $\mathcal{C}$ at $x_{0}$ is identified with the Lie algebra of right invariant vector fields on $G$. Let 
\begin{equation}
m \colon \pi_{1}(M,x_{0}) \longrightarrow \Aut (\mathcal{C}_{x_{0}})
\end{equation}
be the holonomy homomorphism of $\mathcal{C}$ as a locally constant sheaf over $M$. 

For groups $G_{1}$ and $G_{2}$ and a $G_{1}$-action $\rho$ on $G_{2}$, we define the semidirect product $G_{1} \ltimes_{\rho} G_{2}$ of $G_{1}$ and $G_{2}$ with respect to $\rho$ as the group with the underlying set $G_{1} \times G_{2}$ and the product defined by
\begin{equation}
(g_{1},g_{2}) \cdot (g'_{1},g'_{2}) = (g_{1} \cdot g'_{1}, \rho(g'_{1})(g_{2})  \cdot g'_{2})
\end{equation}
for $g_{1}$, $g'_{1}$ in $G_{1}$ and $g_{2}$, $g'_{2}$ in $G_{2}$. For a subgroup $B$ of $\Aut(G_{1})$, we denote the canonical injection $B \longrightarrow \Aut(G_{2})$ by $\iota$.

For a Lie group $G_{1}$, we denote the Lie algebra of the right invariant vector fields on $G_{1}$ by $\Lie(G_{1})^{-}$. For $a$ in $\Aut(G_{1})$, the tangent map of $a \colon G_{1} \longrightarrow G_{1}$ induces an automorphism $a_{*}$ of $\Lie(G_{1})^{-}$. We define a homomorphism $A \colon \Aut(G_{1}) \longrightarrow \Aut(\Lie(G_{1})^{-})$ by $A(a)v = a_{*}v$.

\begin{Proposition}
There exists a homomorphism
\begin{equation}
\Psi \colon \pi_{1}(M,x_{0}) \ltimes_{\rho_{0}} \pi_{1}(F,x_{0}) \longrightarrow \Aut (G,H) \ltimes_{\iota} G
\end{equation}
such that the diagram 
\begin{equation}\label{Equation : Diagram}
\xymatrix{ \pi_{1}(M,x_{0}) \ar[d]_{i} \ar[rr]^{m} & & \Aut(\mathcal{C}_{x_{0}}) \\
\pi_{1}(M,x_{0}) \ltimes_{\rho_{0}} \pi_{1}(F,x_{0}) \ar[r]_>>>>>{\Psi} & \Aut (G,H) \ltimes_{\iota} G \ar[r]_>>>>>{\pr_{1}} & \Aut (G,H) \ar[u]_{A} }
\end{equation}
commutes where $i$ is the injection defined by $i(\gamma) = (\gamma,e)$ for $\gamma$ in $\pi_{1}(M,x_{0})$, $\pr_{1} \colon \Aut (G,H) \ltimes_{\iota} G \longrightarrow \Aut(G,H)$ is the first projection and $A$ is given by the homomorphism $A \colon \Aut (G,H) \longrightarrow \Aut(\Lie(G)^{-})$ with the identification $\Aut(\Lie(G)^{-}) \cong \Aut(\mathcal{C}_{x_{0}})$.
\end{Proposition}

\begin{proof}
We recall the definition of a homomorphism $\Phi \colon \pi_{0}(\Diff(F,x_{0},\mathcal{F}|_{F})) \longrightarrow \Aut G$ defined in Nozawa \cite{Nozawa}. Since the leaves of $(F,\mathcal{F}|_{F})$ are dense, the Lie algebra of transverse fields 
\begin{equation}
C_{b}^{\infty}(TF/T(\mathcal{F}|_{F})) = \{ Y \in C^{\infty}(TF/T(\mathcal{F}|_{F})) \ | \ [Z,Y]=0, \forall Z \in C^{\infty}(T(\mathcal{F}|_{F})) \}
\end{equation}
is equal to the structural Lie algebra $\Lie(G)$ of the Lie foliation $(F,\mathcal{F}|_{F})$ (see Remark in page 117 of \cite{Molino}). For $f$ in $\Diff(F,x_{0},\mathcal{F}|_{F})$ and a transverse field $X$ on $(F,\mathcal{F}|_{F})$, we have $[Y,f_{*}X] = f_{*}[f^{-1}_{*}Y,X]=0$ for $Y$ in $C^{\infty}(T(\mathcal{F}|_{F}))$. Hence $f_{*}X$ is a transverse field. Thus $f_{*}$ induces an automorphism of $\Lie(G)$. By the Lie's theorem, there exists a unique automorphism $f_{G}$ of $G$ which induces $f_{*}$ on $\Lie(G)$. We define $\Phi(f)=f_{G}$. This $\Phi$ is clearly a homomorphism. As shown in Nozawa \cite{Nozawa}, $f_{G}$ depends only on the isotopy class of $f$.

We define $\Psi \colon \pi_{1}(M,x_{0}) \ltimes_{\rho_{0}} \pi_{1}(F,x_{0}) \longrightarrow \Aut (G,H) \ltimes_{\iota} G$ by
\begin{equation}
\Psi(\gamma_{1},\gamma_{2})=(\Phi \circ h_{p}(\gamma_{1}),\hol_{F}(\gamma_{2}))
\end{equation}
for $\gamma_{1}$ in $\pi_{1}(M,x_{0})$ and $\gamma_{2}$ in $\pi_{1}(F,x_{0})$. We show that $\Psi$ is a homomorphism. For $\gamma_{1}$, $\gamma'_{1}$ in $\pi_{1}(M,x_{0})$ and $\gamma_{2}$, $\gamma'_{2}$ in $\pi_{1}(F,x_{0})$, we have
\begin{equation}\label{Equation : Homomorphism1}
\Psi(\gamma_{1},\gamma_{2}) \cdot \Psi(\gamma'_{1},\gamma'_{2}) = (\Phi \circ h_{p}(\gamma_{1}) \cdot \Phi \circ h_{p}(\gamma'_{1}), \iota(\Phi \circ h_{p}(\gamma'_{1}))(\hol_{F}(\gamma_{2})) \cdot \hol_{F}(\gamma'_{2}))
\end{equation}
and
\begin{equation}\label{Equation : Homomorphism2}
\Psi((\gamma_{1},\gamma_{2}) \cdot (\gamma'_{1},\gamma'_{2})) = (\Phi \circ h_{p}(\gamma_{1} \cdot \gamma'_{1}) , \hol_{F} (\rho_{0}(\gamma'_{1})(\gamma_{2}) \cdot \gamma'_{2})).
\end{equation}
Since $\Phi$ and $h_{p}$ are homomorphisms, the first components of the right hand sides of \eqref{Equation : Homomorphism1} and \eqref{Equation : Homomorphism2} are equal. To show that $\Psi$ is a homomorphism, it suffices to show that the second components of the right hand sides of \eqref{Equation : Homomorphism1} and \eqref{Equation : Homomorphism2} are equal. Since $\hol_{F}$ is a homomorphism, it suffices to show
\begin{equation}\label{Equation : Homomorphism3}
\iota(\Phi \circ h_{p}(\gamma'_{1}))(\hol_{F}(\gamma_{2})) = \hol_{F} (\rho_{0}(\gamma'_{1})(\gamma_{2})).
\end{equation}
We put $[f]=h_{p}(\gamma'_{1})$. Since $\rho_{0}=\pi_{1} \circ h_{p}$, \eqref{Equation : Homomorphism3} is equivalent to 
\begin{equation}\label{Equation : Homomorphism4}
\iota(\Phi([f]))(\hol_{F}(\gamma_{2})) = \hol_{F} (\pi_{1}([f])(\gamma_{2})).
\end{equation}
This is equivalent to the Eq. 9 in Nozawa \cite{Nozawa} shown there.

The commutativity of the diagram \eqref{Equation : Diagram} follows from the commutativity of the diagram (7) in Nozawa \cite{Nozawa}. In fact, the commutativity of the diagram implies that
\begin{equation}\label{Equation : m}
m(\gamma)=A \circ \Phi \circ h_{p}(\gamma)
\end{equation}
for every $\gamma$ in $\pi_{1}(M,x_{0})$, which implies the commutativity of the diagram \eqref{Equation : Diagram}.

We show that the image of $\Psi$ is contained in $\Aut(G,H) \ltimes_{\iota} G$. By Lemma \ref{Lemma : Pi 2}, the subgroup $\partial (\pi_{2}(W,p(x_{0})))$ of $\pi_{1}(F,x_{0})$ is preserved by $\rho_{0}$. Hence $\pi_{1}([f])(\gamma)$ is an element of $\partial (\pi_{2}(W,p(x_{0})))$ for $\gamma$ in $\partial (\pi_{2}(W,p(x_{0})))$ and $f$ in $\Diff(F,x_{0},\mathcal{F}|_{F})$. By \eqref{Equation : Homomorphism4}, we have $\iota(\Phi(f))(\hol_{F}(\gamma)) = \hol_{F} (\pi_{1}([f])(\gamma))$. Hence $\iota(\Phi(f))$ preserves $\hol_{F} (\partial (\pi_{2}(W,p(x_{0}))))$.
 Since $H$ is the closure of $\hol_{F} (\partial (\pi_{2}(W,p(x_{0}))))$ in $G$, $\iota(\Phi(f))$ preserves $H$. Thus the image of $\Psi$ is contained in $\Aut(G,H) \ltimes_{\iota} G$.
\end{proof}

Let $\Gamma$ be a group. A linear presentation $f \colon \Gamma \longrightarrow \Aut(V)$ of $\Gamma$ on a vector space $V$ is defined to be unipotent if there exists $n$ such that $(\id_{V} - f(g_{1})) \circ (\id_{V} - f(g_{2})) \circ \cdots \circ (\id_{V} - f(g_{n}))$ is the zero map on $V$ for every $n$-tuple of elements $g_{1}$, $g_{2}$,$\,\ldots\,$, $g_{n}$ of $\Gamma$. 

Let $G_{1}$ be a connected simply connected Lie group. Let $S$ be a subgroup of $\Aut(G_{1})$. Let $I$ be a dense subgroup of $G_{1}$ preserved by the action of $S$. We put $J=S \ltimes_{\iota} I$. We have
\begin{Lemma}\label{Lemma : Nilpotency}
$J$ is nilpotent if and only if $G_{1}$ is nilpotent and $A|_{S} \colon S \longrightarrow \Lie(G_{1})$ is unipotent. 
\end{Lemma}

\begin{proof}
Assume that $J$ is nilpotent. Since the nilpotency of topological groups is a closed condition, $\overline{J}$ is also nilpotent. Obviously $G_{1}$ is nilpotent. By the Engel's theorem and the nilpotency of $\overline{J}$, the adjoint action of $\overline{J}$ on $\Lie(\overline{J})$ is unipotent. From $\Lie(\overline{J})=\Lie(\overline{S}) \ltimes_{\iota} \Lie(G_{1})$, the action of $\overline{S}$ on $\Lie(G_{1})$ is unipotent. Hence the action of $S$ on $\Lie(G_{1})$ is unipotent.

Conversely assume that $G_{1}$ is nilpotent and $A|_{S} \colon S \longrightarrow \Lie(G_{1})$ is unipotent. By Engel's theorem, the adjoint action of $G_{1}$ on $\Lie(G_{1})$ is unipotent. Hence $\Ad g$ is a unipotent element of $\Aut(G_{1})$ for every $g$ in $G_{1}$. By the assumption on the unipotency of $A|_{S}$, $\Ad a$ is a unipotent element of $\Aut(\Lie(G_{1}))$ for every $a$ in $\overline{S}$. Since $(a,g) = (\id_{G_{1}}, g) \circ (a,e)$ in $\overline{J}$ and the product of two unipotent elements is unipotent, $\Ad(a,g)$ is unipotent for every element $(a,g)$ of $\overline{J}$. By Engel's theorem, $\overline{J}$ is nilpotent. Hence $J$ is nilpotent.
\end{proof}

We denote the injection $F \longrightarrow M$ by $i_{F}$. We consider $\Psi_{H}$ and $\Psi_{G/H}$ in the following diagram:
\begin{equation}
\xymatrix{ 
 & \Aut (H) \ltimes_{\iota} H \\
\pi_{1}(M,x_{0}) \ltimes_{\rho_{0}} \partial (\pi_{2}(W,p(x_{0}))) \ar[r] \ar[d] \ar[ur]^{\Psi_{H}} & \Aut (G,H) \ltimes_{\iota} H \ar[d] \ar[u] \\
 \pi_{1}(M,x_{0}) \ltimes_{\rho_{0}} \pi_{1}(F,x_{0}) \ar[r]^>>>>>>>>>>{\Psi} \ar[d]_{\id \times i_{F *}} & \Aut (G,H) \ltimes_{\iota} G \ar[d] \\
 \pi_{1}(M,x_{0}) \ltimes_{\rho_{0}} i_{F *}(\pi_{1}(F,x_{0})) \ar[r] \ar[r]_<<<<<<{\Psi_{G/H}} & \Aut (G/H) \ltimes_{\iota} (G/H) }
\end{equation}
where the vertical arrows are canonical injections or canonical projections. Note that $\partial (\pi_{2}(W,p(x_{0})))$ is preserved by $\rho_{0}$ by Lemma \ref{Lemma : Pi 2}.
\begin{Proposition}\label{Proposition : Characterization}
For a subgroup $K$ of $\pi_{1}(M,x_{0})$ of finite index, $m|_{K}$ is unipotent  and $G/H$ is nilpotent if and only if $\Psi_{H}(K \ltimes_{\rho_{0}} \partial (\pi_{2}(W,p(x_{0}))))$ and $\Psi_{G/H}(K \ltimes_{\rho_{0}} (i_{F *}(\pi_{1}(F,x_{0})) \cap K))$ are nilpotent.
\end{Proposition}

\begin{proof}
We define $m_{H} = \pi_{H} \circ m$ and $m_{G/H} = \pi_{G/H} \circ m$ where $\pi_{H} \colon \Aut(G,H) \longrightarrow \Aut(H)$ and $\pi_{G/H} \colon \Aut(G,H) \longrightarrow \Aut(G/H)$ are the canonical maps. 

An automorphism $a$ of $\Lie(G)^{-}$ which preserves $\Lie(H)^{-}$ is zero if and only if $a|_{\Lie(H)^{-}}$ is zero and $a$ induces the zero map on $\Lie(G/H)^{-}$. Then there exists $n$ such that $(\id_{\mathcal{C}_{x_{0}}} - m(g_{1})) \circ (\id_{\mathcal{C}_{x_{0}}} - m(g_{2})) \circ \cdots \circ (\id_{\mathcal{C}_{x_{0}}} - m(g_{n}))$ is the zero map on $\Lie(\mathcal{C}_{x_{0}})$ for every $n$-tuple $\{g_{i}\}_{i=1}^{n}$ of elements of $K$ if and only if there exists $n$ such that
\begin{enumerate}
\item The map
\begin{equation}
(\id_{\Lie(H)^{-}} - m_{H}(g_{1})) \circ (\id_{\Lie(H)^{-}} - m_{H}(g_{2})) \circ \cdots \circ (\id_{\Lie(H)^{-}} - m_{H}(g_{n}))
\end{equation}
is the zero map on $\Lie(H)^{-}$ and
\item The map
\begin{equation}
(\id_{\Lie(G/H)^{-}} - m_{G/H}(g_{1})) \circ (\id_{\Lie(G/H)^{-}} - m_{G/H}(g_{2})) \circ \cdots \circ (\id_{\Lie(G/H)^{-}} - m_{G/H}(g_{n}))
\end{equation}
induces the zero map on $\Lie(G/H)^{-}$.
\end{enumerate}
Hence $m|_{K}$ is unipotent if and only if $m_{H}|_{K}$ and $m_{G/H}|_{K}$ are unipotent. Since $\hol_{F}(\pi_{1}(F,x_{0}))$ is dense in $G$ by Lemma 4.23 of \cite{Moerdijk Mrcun}, $\hol_{F}(\pi_{1}(F,x_{0}))/\overline{\hol_{F}(K)}$ is dense in $G/\overline{\hol_{F}(K)}$. Since $G$ is connected and $K$ is of finite index in $\pi_{1}(F,x_{0})$, $G/\overline{\hol_{F}(K)}$ is trivial. Hence $\hol_{F}(K)$ is dense in $G$. Then $\hol_{F}(K)/H$ is dense in $G/H$. Hence, if we put $G_{1}= H$, $S =\Psi_{H}(K \ltimes_{\rho_{0}} \{e\})$ and $I = \hol_{F}(\partial (\pi_{2}(W,p(x_{0}))))$ or $G_{1}=G/H$, $S = \Psi_{G/H}(K \ltimes_{\rho_{0}} \{e\})$ and $I = \hol_{F}(K)/H$ in Lemma \ref{Lemma : Nilpotency}, the assumption of Lemma \ref{Lemma : Nilpotency} are satisfied. Then, by Lemma \ref{Lemma : Nilpotency}, $G/H$ is nilpotent and $m_{H}|_{K}$ and $m_{G/H}|_{K}$ are unipotent if and only if $\Psi_{H}(K \ltimes_{\rho_{0}} \partial (\pi_{2}(W,p(x_{0}))))$ and $\Psi_{G/H}(K \ltimes_{\rho_{0}} (i_{F *}(\pi_{1}(F,x_{0})) \cap K))$ are nilpotent.
\end{proof}

Let $K$ be a subgroup of $\pi_{1}(M,x_{0})$. Among $\Psi_{H}(K \ltimes_{\rho_{0}} \partial (\pi_{2}(W,p(x_{0}))))$ and $\Psi_{G/H}(K \ltimes_{\rho_{0}} (i_{F *}(\pi_{1}(F,x_{0})) \cap K))$, the nilpotency of the latter follows from the nilpotentcy of $K$ as shown in the following lemma. On the other hand, the nilpotency of $\Psi_{H}(K \ltimes_{\rho_{0}} \partial (\pi_{2}(W,p(x_{0}))))$ is not controlled by the nilpotency of $\pi_{1}(M,x_{0})$ as shown in Example \ref{Example}.
\begin{Lemma}\label{Lemma : Nilpotency of Semidirectproduct}
If $K$ is nilpotent, then $K \ltimes_{\rho_{0}} (i_{F *}(\pi_{1}(F,x_{0})) \cap K)$ is nilpotent.
\end{Lemma}

\begin{proof}
For $\delta$ in $\pi_{1}(F,x_{0})$ and $\gamma$ in $\pi_{1}(M,x_{0})$, we have $i_{F *}[h_{p}(\gamma)\delta]=i_{F *}[\gamma \cdot \delta \cdot \gamma^{-1}]$ in $\pi_{1}(M,x_{0})$. Since $\rho_{0}=\pi_{1} \circ h_{p}$, the action $\rho_{0}|_{K} \colon K \longrightarrow \Aut(i_{F *}(\pi_{1}(F,x_{0})) \cap K)$ induced by $\rho_{0}$ is the action by conjugation. 

Let $\Gamma$ be a nilpotent group and $\Pi$ be a normal subgroup of $\Gamma$. Let $\rho_{1} \colon \Gamma \longrightarrow \Aut(\Pi)$ be the action of $\Gamma$ on $\Pi$ defined by $\rho_{1}(g)(h)=g \cdot h \cdot g^{-1}$ for $g$ in $\Gamma$ and $h$ in $\Pi$. To show Lemma \ref{Lemma : Nilpotency of Semidirectproduct}, it suffices to show that $\Gamma \ltimes_{\rho_{1}} \Pi$ is nilpotent. Put $\tilde{\Gamma} = \Gamma \ltimes_{\rho_{1}} \Pi$. For a group $C$, write $C_{0}=C$ and $C_{i} = [C,C_{i-1}]$. Let $\pr_{1} \colon \tilde{\Gamma} \longrightarrow \Gamma$ be the first projection. Since $\pr_{1} (\tilde{\Gamma}_{i})=\Gamma_{i}$ and by the nilpotency of $\Gamma$, there exists $i_{0}$ such that $\tilde{\Gamma}_{i_{0}}$ is contained in $\{e\} \ltimes_{\rho_{1}} \Pi$. For $(g_{1},h_{1})$ and $(e,h_{2})$ in $\tilde{\Gamma}$, we have
\begin{equation}
[(g_{1},h_{1}),(e,h_{2})] = (e,[g_{1}^{-1}h_{1},h_{2}]).
\end{equation}
Hence $[\tilde{\Gamma},\{e\} \ltimes_{\rho_{1}} (\Pi \cap \Gamma_{i})]$ is contained in $\{e\} \ltimes_{\rho_{1}} (\Pi \cap \Gamma_{i+1})$. Hence there exists $i_{1}$ such that $\tilde{\Gamma}_{i_{1}}=\{e\}$ by the nilpotency of $\Gamma$.
\end{proof}

The developability of $(M,\mathcal{F})$ implies the triviality of  $\Psi_{H}(K \ltimes_{\rho_{0}} \partial (\pi_{2}(W,p(x_{0}))))$ according to the following lemma:  
\begin{Lemma}[Molino, Exercise 4.7 {[4]} of \cite{Molino}]\label{Lemma : Developability}
$(M,\mathcal{F})$ is developable if and only if the image of $\hol_{F} \circ \partial \colon \pi_{2}(W,p(x_{0})) \longrightarrow G$ is discrete.
\end{Lemma}

\begin{proof}
Let $\pi \colon \tilde{M} \longrightarrow M$ be the universal cover of $M$. Let $\tilde{F}$ be a connected component of the fiber of $F$. Since $(\tilde{F},(\pi^{*}\mathcal{F})|_{\tilde{F}})$ is a covering space of $(F,\mathcal{F}|_{F})$, $\pi^{*}\mathcal{F}|_{\tilde{F}}$ is a $G$-Lie foliation on $\tilde{F}$. Fix a point $\tilde{x}_{0}$ on the fiber of $x_{0}$. The holonomy homomorphism of the Lie foliation $(\tilde{F},(\pi^{*}\mathcal{F})|_{\tilde{F}})$ is given by 
\begin{equation}
\xymatrix{
\pi_{1}(\tilde{F},\tilde{x}_{0}) \ar[r]^{\pi_{*}} & \pi_{1}(F,x_{0}) \ar[r]^<<<<<{\hol_{F}} & G. 
}
\end{equation}
By the definition of $\pi$ and the homotopy exact sequence, we have $\pi_{*}\pi_{1}(\tilde{F},\tilde{x}_{0}) = \ker i_{F *} = \partial (\pi_{2}(W,p(x_{0})))$. Let $H$ be the closure of $\hol_{F}(\pi_{*}(\pi_{1}(\tilde{F},\tilde{x}_{0})))$ in $G$. Let $\overline{L}$ be the closure of a leaf $L$ of $\pi^{*}\mathcal{F}$. Then $(\overline{L},(\pi^{*}\mathcal{F})|_{\overline{L}})$ is an $H$-Lie foliation with dense leaves by Theorem 4.2 of Molino \cite{Molino}.

Clearly the image of $\hol_{F} \circ \partial$ is discrete in $G$ if and only if $H$ is of dimension $0$. Since $\dim H$ is equal to the codimension of $(\overline{L},(\pi^{*}\mathcal{F})|_{\overline{L}})$, $H$ is of dimension $0$ if and only if $(\overline{L},(\pi^{*}\mathcal{F})|_{\overline{L}})$ is a foliation with one leaf $\overline{L}$. Since $\overline{L}$ is the closure of a leaf $L$ of $\pi^{*}\mathcal{F}$, $(\overline{L},(\pi^{*}\mathcal{F})|_{\overline{L}})$ is a foliation with one leaf $\overline{L}$ if and only if $L$ is closed in $\tilde{M}$. By the homogeneity of the transversely complete foliation $(\tilde{M},\pi^{*}\mathcal{F})$ (see Theorem 4.2' of Molino \cite{Molino}), a leaf $L$ of $(\tilde{M},\pi^{*}\mathcal{F})$ is closed in $\tilde{M}$ if and only if every leaf of $\pi^{*}\mathcal{F}$ is closed in $\tilde{M}$. Since the closure of each leaf of $\pi^{*}\mathcal{F}$ is a fiber of the basic fibration of $(\tilde{M},\pi^{*}\mathcal{F})$ by Theorem 4.2' of \cite{Molino}, every leaf of $\pi^{*}\mathcal{F}$ is closed in $\tilde{M}$ if and only if $(\tilde{M},\pi^{*}\mathcal{F})$ is developable.
\end{proof}

We obtain a sufficient condition for the unipotency of the holonomy homomorphism of the Molino's commuting sheaves of developable Riemannian foliations. 
\begin{Proposition}\label{Proposition : Unipotency}
Let $(M,\mathcal{F})$ be a developable Riemannian foliation on a closed manifold. Let $K$ be a subgroup of of $\pi_{1}(M,x_{0})$ of finite index. If $K$ is nilpotent, then the restriction of the holonomy homomorphism of the Molino's commuting sheaf $\mathcal{C}$ of $(M,\mathcal{F})$ to $K$ is unipotent.
\end{Proposition}

\begin{proof}
First, we show Proposition \ref{Proposition : Unipotency} in the case where $(M,\mathcal{F})$ is transversely parallelizable. Let $\mathcal{F}$ be a developable transversely parallelizable foliation on $M$. Let $p \colon M \longrightarrow W$ be the basic fiberation of $(M,\mathcal{F})$. Fix a fiber $F$ of $p$ and a point $x_{0}$ on $F$. Let $m \colon \pi_{1}(M,x_{0}) \longrightarrow \Aut(\mathcal{C}_{x_{0}})$ be the holonomy homomorphism of $\mathcal{C}$. Let $\partial \colon \pi_{2}(F,p(x_{0})) \longrightarrow \pi_{1}(M,x_{0})$ be the connecting homomorphism in the homotopy exact sequence of the fiber bundle $p$. By the developability of $(M,\mathcal{F})$ and Lemma \ref{Lemma : Developability}, the image of $\xymatrix{ \pi_{2}(W,p(x_{0})) \ar[r]^<<<<<{\partial} & \pi_{1}(F,x_{0}) \ar[r]^<<<<<{\hol_{F}} & G }$ is discrete. Hence we have $H=\{e\}$. Hence $\Psi_{H}(K \ltimes_{\rho_{0}} \partial \pi_{2}(W,p(x_{0})))$ is trivial. Since $K$ is nilpotent, $K \ltimes_{\rho_{0}} (i_{F *}(\pi_{1}(F,x_{0})) \cap K)$ is nilpotent by Lemma \ref{Lemma : Nilpotency of Semidirectproduct}. Hence $\Psi_{G/H}(K \ltimes_{\rho_{0}} (i_{F *}(\pi_{1}(F,x_{0})) \cap K))$ is nilpotent. Then $m|_{K}$ is unipotent by Proposition \ref{Proposition : Characterization}. 

We consider the general case. Let $\mathcal{F}$ be a developable Riemannian foliation on $M$. Let $p \colon M^{1} \longrightarrow M$ be the transverse orthonormal frame bundle of $(M,\mathcal{F})$. Let $\mathcal{F}^{1}$ be the lift of $\mathcal{F}$ to $M^{1}$. Let $\mathcal{C}$ and $\mathcal{C}^{1}$ be the Molino's commuting sheaves of $(M,\mathcal{F})$ and $(M^{1},\mathcal{F}^{1})$ respectively. Recall that, by the definition of the Molino's commuting sheaf of Riemannian foliations (see Section 5.3 of \cite{Molino}), the stalk $\mathcal{C}_{x}$ at $x$ on $M$ is the set of the germs of vector fields at $x$ whose canonical lift to $M^{1}$ is a section of $\mathcal{C}^{1}$. Then we have $p^{*}\mathcal{C} = \mathcal{C}^{1}$. Let $m$ and $m^1$ be the holonomy homomorphisms of $\mathcal{C}$ and $\mathcal{C}^{1}$ respectively. Since $p^{*}\mathcal{C} = \mathcal{C}^{1}$, $m^{1}$ is given by
\begin{equation}
\xymatrix{ \pi_{1}(M^{1},\tilde{x}_{0}) \ar[r]^{p_{*}} & \pi_{1}(M,x_{0}) \ar[r]^{m} & \Aut (\mathcal{C}_{x_{0}}) \ar[r] & \Aut (\mathcal{C}^{1}_{\tilde{x}_{0}}) } 
\end{equation}
where the last map is the identification by the pullback by $p$. Note that $p_{*}^{-1}(K)$ is nilpotent. By the previous case, $m^{1}|_{p^{-1}(K)}$ is unipotent. Hence $m|_{K}$ is unipotent.
\end{proof}

\section{Vanishing of secondary characteristic classes of flat vector bundles with unipotent holonomy homomorphisms}

We use the following computation of secondary characteristic classes of flat vector bundles by Kamber and Tondeur to show the nullity of secondary characteristic classes of flat vector bundles with unipotent holonomy homomorphisms.
\begin{Theorem}[Kamber-Tondeur, Eqs. 5.74 and 6.31 in \cite{Kamber Tondeur}]
\begin{equation}\label{Equation : Cohomology}
H^{\bullet}(\mathfrak{gl}_{r},O_{r}) \cong \wedge^{\bullet} (y_{1},y_{3},\ldots,y_{r'}),
\end{equation}
where $r'$ is the largest odd integer equal to or less than $r$, and $y_{2i-1}$ is of degree $2i-1$ and represented by the cocycle
\begin{equation}\label{Equation : Cocycle}
\left( -\frac{1}{2} \right)^{i-1} \frac{(i-1)! i!}{(2i-1)!} c_{i}([\theta,\theta] \wedge \cdots \wedge [\theta,\theta] \wedge \theta).
\end{equation}
Here, $c_{i}$ is the $i$-th Chern polynomial regarded as an element of $S^{i} \mathfrak{gl}_{r}^{*}$, and $\theta$ is the identity map on $\mathfrak{gl}_{r}$, $[\theta,\theta]$ is defined by $[\theta,\theta](x,y)=[x,y]$ for $x$, $y$ in $\mathfrak{gl}_{r}$, and $[\theta,\theta] \wedge \cdots \wedge [\theta,\theta] \wedge \theta$ is regarded as an element of $\Hom(\wedge^{2i-1} \mathfrak{gl}_{r}, S^{i} \mathfrak{gl}_{r})$.
\end{Theorem}

Let $M$ be a closed manifold. Let $E$ be a vector bundle over $M$ of rank $r$ with a flat connection $\omega$. The secondary characteristic classes of $(E,\omega)$ are the cohomology classes in the image of the generalized characteristic homomorphism
\begin{equation}
\Delta \colon H^{\bullet}(\mathfrak{gl}_{r},\Or_{r}) \longrightarrow H^{\bullet}(M;\mathbb{R})
\end{equation}
of $(E,\Delta)$ defined as follows (see \cite{Kamber Tondeur}): Let $P$ be the frame bundle of $E$. By the flat connection form on $P$ associated to $\omega$, we have a linear map
\begin{equation}
\omega \colon \mathfrak{gl}_{r}^{*} \longrightarrow C^{\infty}(T^{*}P). 
\end{equation}
This $\omega$ induces
\begin{equation}
\omega \colon \wedge^{\bullet} \mathfrak{gl}_{r}^{*} \longrightarrow C^{\infty}(\wedge^{\bullet} T^{*}P). 
\end{equation}
Let $P/\Or_{r}$ be the $(\GL_{r}/\Or_{r})$-bundle over $M$ obtained by taking quotient of $P$ by the $\Or_{r}$-subaction of the principal $\GL_{r}$-action. Fix an $\Or_{r}$-reduction $\Or$ of $P$. This $\Or$ determines a section $s \colon M \longrightarrow P/\Or_{r}$. Let $(\wedge^{\bullet} (\mathfrak{gl}_{r}^{*}/\mathfrak{o}_{r}^{*}))^{\Or_{r}}$ be the subspace of $\Or_{r}$-invariant elements of $\wedge^{\bullet} (\mathfrak{gl}_{r}^{*}/\mathfrak{o}_{r}^{*})$. Then $\omega$ induces a map 
\begin{equation}\label{Equation : Omega}
\omega \colon (\wedge^{\bullet} (\mathfrak{gl}_{r}^{*}/\mathfrak{o}_{r}^{*}))^{\Or_{r}} \longrightarrow C^{\infty}( \wedge^{\bullet} T^{*}(P/\Or_{r})) 
\end{equation}
By the definition, $H^{\bullet}(\mathfrak{gl}_{r},\Or_{r})$ is the cohomology of the chain complex $(\wedge^{\bullet} (\mathfrak{gl}_{r}^{*}/\mathfrak{o}_{r}^{*}))^{\Or_{r}}$ with the differential induced by the Eilenberg-MacLane differential. Then the composition of \eqref{Equation : Omega} and the pullback by $s$
\begin{equation}
s^{*} \colon C^{\infty}( \wedge^{\bullet} T^{*}(P/\Or_{r})) \longrightarrow C^{\infty}( \wedge^{\bullet} T^{*}M)
\end{equation}
induces $\Delta \colon H^{\bullet}(\mathfrak{gl}_{r},\Or_{r}) \longrightarrow H^{\bullet}(M;\mathbb{R})$. Note that the contractibility of the fibers of $P/\Or_{r}$ implies that the secondary characteristic classes of $(E,\omega)$ are independent of $\Or$.

We show
\begin{Proposition}\label{Proposition : Vanishing}
Fix a point $x_{0}$ on $M$. Let $K$ be a subgroup of $\pi_{1}(M,x_{0})$ of finite index. If the restriction of the holonomy homomorphism 
\begin{equation}
h_{E} \colon \pi_{1}(M,x_{0}) \longrightarrow \Aut E_{x_{0}}
\end{equation}
of $E$ to $K$ is unipotent, then the secondary characteristic classes of $(E,\omega)$ are zero in $H^{\bullet}(M;\mathbb{R})$.
\end{Proposition}

\begin{proof}
Let $\pi \colon M' \longrightarrow M$ be the covering space of $M$ such that $\pi_{*} \pi_{1}(M',x'_{0}) = K$. Since $K$ is of finite index, $\pi$ is a finite covering. Hence $M'$ is also closed. By the naturality of the secondary characteristic classes, the secondary characteristic classes of $(\pi^{*}E,\pi^{*}\omega)$ are the pullback of the secondary characteristic classes of $(E, \omega)$ by $\pi^{*}$. Since $\pi^{*} \colon H^{\bullet}(M ; \mathbb{R}) \longrightarrow H^{\bullet}(M' ; \mathbb{R})$ is injective, to show that the secondary characteristic classes of $(E, \omega)$ vanish, it suffices to show that the secondary characteristic classes of $(\pi^{*}E,\pi^{*}\omega)$ vanish. Hence it suffices to show the case of $K=\pi_{1}(M,x_{0})$ for the proof of Proposition \ref{Proposition : Vanishing}.

We assume $K=\pi_{1}(M,x_{0})$. Let $H$ be a connected Lie subgroup of $\Aut (E_{x_{0}})$ which is invariant by the action $\pi_{1}(M,x_{0}) \longrightarrow \Aut(E_{x_{0}})$. Fix a frame $f_{0}$ of $E_{x_{0}}$. By the parallel transport of $Hf_{0}$ with respect to the flat connection, we obtain an $H$-reduction $P_{H}$ of $E$ with a flat connection.

Let 
\begin{equation}
\omega_{H} \colon \mathfrak{h}^{*} \longrightarrow C^{\infty}(T^{*}P_{H}). 
\end{equation}
be the map determined by the flat $H$-connection on $P_{H}$. Then we have a commutative diagram
\begin{equation}\label{Equation : gl and h}
\xymatrix{ \wedge^{\bullet}\mathfrak{gl}_{r}^{*} \ar[r]^>>>>>>{\omega} \ar[d]_{j} & C^{\infty}(\wedge^{\bullet} T^{*}P) \ar[d]^{i} \\
 \wedge^{\bullet} \mathfrak{h}^{*} \ar[r]_>>>>>{\omega_{H}} & C^{\infty}(\wedge^{\bullet} T^{*}P_{H}) }
\end{equation}
where $j$ is induced by the restriction map $\mathfrak{gl}_{r}^{*} \longrightarrow \mathfrak{h}^{*}$ and $i$ is the restriction map.

The diagram \eqref{Equation : gl and h} induces a commutative diagram
\begin{equation}
\xymatrix{ (\wedge^{\bullet}(\mathfrak{gl}_{r}^{*}/\mathfrak{o}_{r}^{*}))^{\Or_{r}} \ar[r]^>>>>>>>>>>>{\omega} \ar[d]_{j} & C^{\infty}(\wedge^{\bullet} T^{*}(P/\Or_{r})) \ar[d]^{i} \\
 (\wedge^{\bullet} (\mathfrak{h}^{*}/(\mathfrak{h} \cap \mathfrak{o}_{r})^{*}))^{H \cap \Or_{r}} \ar[r]_{\omega_{H}} & C^{\infty}(\wedge^{\bullet} T^{*}(P_{H}/(\Or_{r} \cap H)))}
\end{equation}
where $(\wedge^{\bullet} \mathfrak{h}_{r}^{*}/(\mathfrak{h} \cap \mathfrak{o}_{r})^{*})^{H \cap \Or_{r}}$ is the subspace of $(H \cap \Or_{r})$-invariant elements of $\wedge^{\bullet} (\mathfrak{h}_{r}^{*}/(\mathfrak{h} \cap \mathfrak{o}_{r})^{*})$. 

Let $N_{r}$ be the Lie subgroup of $\GL_{r}$ consisting of unipotent matrices. Then $N_{r}$ is contractible and $\Lie(N_{r}) \cap \mathfrak{o}_{r}=\{0\}$. We have
\begin{equation}
\xymatrix{ (\wedge^{\bullet}(\mathfrak{gl}_{r}^{*}/\mathfrak{o}_{r}^{*}))^{\Or_{r}} \ar[r]^>>>>>{\omega} \ar[d]_{j} & C^{\infty}(\wedge^{\bullet} T^{*}(P/\Or_{r})) \ar[d]^{i} \\
 \wedge^{\bullet} \Lie(N_{r})^{*} \ar[r]_>>>>>>{\omega_{H}} & C^{\infty}(\wedge^{\bullet} T^{*}(P_{N_{r}}))}
\end{equation}
The representative of the generators $y_{1}$, $y_{2}$,$\,\ldots\,$, $y_{r'}$ of $H^{\bullet}(\mathfrak{gl}_{r},\Or_{r})$ are mapped to zero by $j \colon \wedge^{\bullet}(\mathfrak{gl}_{r}^{*}/\mathfrak{o}_{r}^{*}) \longrightarrow \wedge^{\bullet} \Lie(N_{r})^{*}$. In fact, since Chern polynomials are functions of eigenvalues of the elements of Lie algebras as matrices, their restrictions to $\wedge^{\bullet} \Lie(N_{r})^{*}$ vanish. $\pi(y_{i})$ is zero as a cocycle in $\wedge^{\bullet} \Lie(N_{r})^{*}$. Since $i$ is a homotopy equivalence by the contractibility of $N_{r}$, $i$ induces an isomorphism on the cohomology. Hence Proposition \ref{Proposition : Vanishing} is proved.
\end{proof}

\section{Proof of Theorems and Corollaries}

We show Theorem \ref{Theorem : Vanishing} by using a theorem of Gromov \cite{Gromov}.
\begin{proof}[Proof of Theorem \ref{Theorem : Vanishing}]
By a theorem of Gromov \cite{Gromov}, there exists a nilpotent subgroup $K$ of $\pi_{1}(M,x_{0})$ of finite index. By Proposition \ref{Proposition : Unipotency}, the restriction of the holonomy homomorphism of $\mathcal{C}$ to $K$ is unipotent. Then the secondary characteristic classes of $\mathcal{C}$ vanish by Proposition \ref{Proposition : Vanishing}.
\end{proof}

We show Corollary \ref{Corollary : Minimizability} by using Theorem \ref{Theorem : Vanishing} and theorems of \'{A}lvarez L\'{o}pez. 
\begin{proof}[Proof of Corollary \ref{Corollary : Minimizability}]
Let $\mathcal{C}$ be the Molino's commuting sheaf of $(M,\mathcal{F})$. By Theorem \ref{Theorem : Vanishing}, the secondary characteristic classes of $\mathcal{C}$ vanish. By a theorem of \'{A}lvarez L\'{o}pez \cite{Alvarez Lopez 2}, the \'{A}lvarez class of $(M,\mathcal{F})$ is equal to the secondary characteristic class of degree $1$ up to multiplication of a non-zero constant. By a theorem of \'{A}lvarez L\'{o}pez \cite{Alvarez Lopez}, the \'{A}lvarez class of $(M,\mathcal{F})$ vanishes if and only if $(M,\mathcal{F})$ is minimizable. Hence $(M,\mathcal{F})$ is minimizable.
\end{proof}

By using Lemma \ref{Lemma : Developability}, we show the following proposition.
\begin{Proposition}\label{Lemma : Codimension 3}
Every transversely parallelizable foliation $\mathcal{F}$ of codimension $q$ on a closed manifold $M$ is developable or compact if $q \leq 3$. 
\end{Proposition}

\begin{proof}
Let $F$ be the closure of a leaf of $\mathcal{F}$. By the Molino's structure theorem, $(F,\mathcal{F}|_{F})$ is a $G$-Lie foliation for a Lie group $G$. By the assumption on the codimension of $(M,\mathcal{F})$, we have $\dim M - \dim F \leq 3$.

Assume $\dim F=\dim M$ and $M$ is connected. Since $F$ is a fiber of a basic fibration, we have $M = F$. Then $(M,\mathcal{F})$ is a $G$-Lie foliation. Since every Lie foliation is developable, $(M,\mathcal{F})$ is developable. 

Assume $\dim F=\dim M-1$. Then the base space of the basic fibration of $(M,\mathcal{F})$ is $S^1$. Since $\pi_{2}S^1 = \{1\}$, the image of $\hol_{F} \circ \partial \colon \pi_{2}(S^1,p(x_{0})) \longrightarrow \pi_{1}(F,x_{0}) \longrightarrow G$ is trivial. By Lemma \ref{Lemma : Developability}, $(M,\mathcal{F})$ is developable.

Assume $\dim F=\dim M-2$. Then the base space of the basic fibration of $(M,\mathcal{F})$ is a closed surface $W$. Then $\pi_{2}W$ is a free abelian group of rank less than $2$. Hence the image of $\hol_{F} \circ \partial \colon \pi_{2}(S^1,p(x_{0})) \longrightarrow \pi_{1}(F,x_{0}) \longrightarrow G$ is discrete. By Lemma \ref{Lemma : Developability}, $(M,\mathcal{F})$ is developable.

Assume $\dim F=\dim M-3$. By the assumption on the codimension of $(M,\mathcal{F})$, this implies that $F$ is a leaf of $\mathcal{F}$. Hence $F$ is compact. Since each leaf of $(M,\mathcal{F})$ is a fiber of the basic fibration of $(M,\mathcal{F})$, every leaf of $\mathcal{F}$ is compact.
\end{proof}

The following proposition is well known and follows because Molino's commuting sheaf is the zero sheaf when the leaves are compact \cite{Molino}:

\begin{Proposition}\label{Proposition : Compact}
A compact Riemannian foliation $\mathcal{F}$ on a closed manifold $M$ is minimizable.
\end{Proposition}

We show Theorem \ref{Theorem : Codimension 3} by using Propositions \ref{Lemma : Codimension 3}, \ref{Proposition : Compact} and Corollary \ref{Corollary : Minimizability}.
\begin{proof}[Proof of Theorem \ref{Theorem : Codimension 3}]
Let $\mathcal{F}$ be a transversely parallelizable foliation of codimension $q$ on a closed manifold $M$ with a fundamental group of polynomial growth. Assume that $q \leq 3$. By Proposition \ref{Lemma : Codimension 3}, $(M,\mathcal{F})$ is developable or compact. If $(M,\mathcal{F})$ is developable, then $(M,\mathcal{F})$ is minimizable by Corollary \ref{Corollary : Minimizability}. If $(M,\mathcal{F})$ is compact, then $(M,\mathcal{F})$ is minimizable by Proposition \ref{Proposition : Compact}.
\end{proof}

We deduce Corollary \ref{Corollary : Codimension 2} from Theorem \ref{Theorem : Codimension 3}.
\begin{proof}[Proof of Corollary \ref{Corollary : Codimension 2}]
Let $(M,\mathcal{F})$ be a Riemannian foliation of codimension $2$ on a closed manifold $M$ with a fundamental group of polynomial growth. Let $M^{1}$ be the transverse orthonormal frame bundle of $(M,\mathcal{F})$. Let $\mathcal{F}^{1}$ be the lift of $\mathcal{F}$ to $M^{1}$. Since $\dim M^{1}=\dim M + \dim \Or (2)$ and $\dim \mathcal{F}^{1}= \dim \mathcal{F}$, $\mathcal{F}^{1}$ is of codimesion $3$. Since $\pi_{1}M^{1}$ is a central $\mathbb{Z}$-extension of $\pi_{1}M$ or a central $\mathbb{Z}$-extension of a $\mathbb{Z}/2\mathbb{Z}$-extension of $\pi_{1}M$ depending on the transverse orientability of $(M,\mathcal{F})$, $\pi_{1}M^{1}$ is of polynomial growth. Then, by Theorem \ref{Theorem : Codimension 3}, $(M^{1},\mathcal{F}^{1})$ is minimizable. By Lemma 4.2 of Ghys \cite{Ghys}, $(M,\mathcal{F})$ is minimizable if and only if $(M^{1},\mathcal{F}^{1})$ is minimizable. Hence $(M,\mathcal{F})$ is minimizable.
\end{proof}

\section{Examples}

\subsection{Carri\`{e}re's example}\label{Subsection : Carriere}
Let $A$ be a hyperbolic element of $\SL(2;\mathbb{Z})$. $A$ induces a diffeomorphism $\overline{A}$ of $T^{2}=\mathbb{R}^{2}/\mathbb{Z}^{2}$. Carri\`{e}re's example is the Riemannian flow $\mathcal{F}$ on the mapping torus $T_{A}^{3}=([0,1] \times T^{2})/(0,\overline{A}x) \sim (1,x)$ of $\overline{A}$ whose restriction to each $T^2$-fiber of the fibration $T^2 \longrightarrow T^{3}_{A} \longrightarrow S^1$ is the linear flow tangent to an eigenvector of $A$. These are examples of non-minimizable developable Riemannian flows of codimension $2$ on closed $3$-dimensional manifolds with fundamental groups of exponential growth (see \cite{Carriere}).

\subsection{An example of non-minimizable non-developable Riemannian flow on a closed manifold with abelian fundamental group}\label{Example}
We present an example of a non-minimizable non-developable Riemannian flow on a closed manifold with abelian fundamental group.

Let $(M,\mathcal{F})$ be a closed manifold with a transversely parallelizable foliation. Let $F \longrightarrow M \longrightarrow W$ be the basic fibration of $(M,\mathcal{F})$. Assume that $(M,\mathcal{F})$ is not minimizable and $\pi_{1}M$ is of polynomial growth. There exists a nilpotent subgroup $K$ of $\pi_{1}(M,x_{0})$ of finite index by a theorem of Gromov \cite{Gromov}. We use the notation in the proof of Proposition \ref{Proposition : Characterization}. Since $K$ is nilpotent, $\Psi_{G/H}(K \ltimes_{\rho_{0}} (i_{F *}(\pi_{1}(F,x_{0})) \cap K))$ is nilpotent by Lemma \ref{Lemma : Nilpotency of Semidirectproduct}. If $\Psi_{H}(K \ltimes_{\rho_{0}} \partial (\pi_{2}(W,p(x_{0}))))$ is nilpotent, then $(M,\mathcal{F})$ is minimizable by Propositions \ref{Proposition : Characterization}, \ref{Proposition : Vanishing} and theorems of \'{A}lvarez L\'{o}pez in \cite{Alvarez Lopez} and \cite{Alvarez Lopez 2} as in the proof of Theorem \ref{Theorem : Vanishing} and Corollary \ref{Corollary : Minimizability}. Since $(M,\mathcal{F})$ is not minimizable by the assumption, $\Psi_{H}(K \ltimes_{\rho_{0}} \partial (\pi_{2}(W,p(x_{0}))))$ is not nilpotent by contradiction. By Lemma \ref{Lemma : Nilpotency} in the case where $G_{1}= H$, $S =\Psi_{H}(K \ltimes_{\rho_{0}} \{e\})$ and $I = \hol_{F}(\partial (\pi_{2}(W,p(x_{0}))))$, it follows that $\Psi_{H}(K \ltimes_{\rho_{0}} \partial (\pi_{2}(W,p(x_{0}))))$ is nilpotent if and only if the action $m_{H}|_{K}$ is unipotent. Then $m_{H}|_{K}$ is not unipotent. Since $H$ is the closure of $\hol_{F} \circ \partial(\pi_{2}(W,p(x_{0})))$ in $G$, the action $(\rho_{0}|_{K}) \otimes_{\mathbb{Z}} \mathbb{R}$ of $K$ on $\partial (\pi_{2}(W,p(x_{0}))) \otimes_{\mathbb{Z}} \mathbb{R}$ is not unipotent by \eqref{Equation : Homomorphism3} and \eqref{Equation : m}. Note that $\rho_{0}|_{K}$ is the action by conjugation as shown in the first paragraph of the proof of Lemma \ref{Lemma : Nilpotency of Semidirectproduct}. Thus, if $(M,\mathcal{F})$ is not minimizable and $\pi_{1}M$ is of polynomial growth, then the topology of $M$ is complicated in some sense. The complexity of the action of the diffeomorphism group of the $K3$ surface on its cohomology is relevant in our construction.

Let $\sigma$ be the involution on the $4$-torus $T^{4}=\mathbb{R}^{4}/\mathbb{Z}^{4}$ defined by $\sigma (x_{1},x_{2},x_{3},x_{4}) = (-x_{1},-x_{2},-x_{3},-x_{4})$. Then the quotient space $T^{4}/\sigma$ has $16$ singular points. Blowing up all of the singular points, we have a $K3$ surface $X$.

We will define a $1$-dimensional Riemannian foliation on a total space of a $T^2$-bundle over $X$. Let $\{[dx_{1}], [dx_{2}], [dx_{3}], [dx_{4}]\}$ be the standard basis of $H^{1}(T^{4};\mathbb{Z}) \cong \mathbb{Z}^{4}$. Let $E_{1}$, $E_{2}$, $\cdots$, $E_{16}$ be exceptinal divisors on $X$. We denote the Poincar\'{e} dual of $E_{j}$ by $[E_{j}]$. Since we have $H^{2}(T^{4}/\sigma;\mathbb{Z}) \cong H^{2}(T^{4};\mathbb{Z})^{\sigma} \cong H^{2}(T^{4};\mathbb{Z}) \cong \mathbb{Z}^{6}$, $H^{2}(X;\mathbb{R}) \cong H^{2}(T^4/\sigma;\mathbb{Z}) \oplus \bigoplus_{j=1}^{16}\mathbb{R}[E_{j}] \cong \oplus_{1 \leq k < l \leq 4} \mathbb{R} [dx_{k} \wedge dx_{l}] \oplus \bigoplus_{j=1}^{16}\mathbb{R}[E_{j}]$. For a vector $v=\left( \begin{smallmatrix} v_{1} \\ v_{2} \end{smallmatrix} \right)$ in $\mathbb{Z}^{2}$, let $S_{v}$ be the principal $S^1$-bundle over $X$ whose Euler class is $v_{1}[dx_{1} \wedge dx_{3}] + v_{2}[dx_{2} \wedge dx_{3}]$. Let $\pr_{j} \colon X \times X \longrightarrow X$ be the $j$-th projection and $\delta \colon X \longrightarrow X \times X$ be the diagonal injection. We define a principal $T^2$-bundle $Y$ over $X$ by $Y = \delta^{*} (\pr_{1}^{*}S_{\left( \begin{smallmatrix} 1 \\ 0 \end{smallmatrix} \right)} \oplus \pr_{2}^{*}S_{\left( \begin{smallmatrix} 0 \\ 1 \end{smallmatrix} \right)})$. Let $\Xi \colon \Lie(T^2)=\mathbb{R}^{2} \longrightarrow C^{\infty}(TY)$ be the infinitesimal action of the principal $T^{2}$-action. For $v$ in $\mathbb{R}^{2}$, $\rho_{v}$ denotes the $S^1$-action or the flow on $Y$ generated by $\Xi \left( v \right)$. We denote the $1$-dimensional foliation on $Y$ defined by the orbits of $\rho_{\left( \begin{smallmatrix} 1 \\ \frac{- 1 - \sqrt{5}}{2} \end{smallmatrix}\right)}$ by $\mathcal{G}$. Since $\rho_{\left( \begin{smallmatrix} 1 \\ \frac{- 1 - \sqrt{5}}{2} \end{smallmatrix}\right)}$ preserves any $T^2$-invariant metric on $Y$, $\mathcal{G}$ is a Riemannian foliation.

We will construct a diffeomorphism on $Y$ preserving $\mathcal{G}$. For a primitive vector $v$ in $\mathbb{Z}^{2}$, let $G_{v}$ be the isotropy subgroup of $T^{2}$ corresponding to $v$, which is isomorphic to $S^{1}$. We identify $T^{2}/G_{v}$ with $S^{1}$ by choosing a vector $w$ in $\mathbb{R}^{2}$ so that $\det \left( \begin{matrix} v & w \end{matrix}\right) > 0$ and $w$ is a generator of $S^{1}$. By this identification, for $v$ in $\mathbb{Z}^{2}$, $Y/\rho_{v}$ is considered as a principal $S^1$-bundle over $X$. We show that $Y/\rho_{\left( \begin{smallmatrix} v_{1} \\ v_{2} \end{smallmatrix} \right)}$ is isomorphic to $S_{\left( \begin{smallmatrix} v_{1} \\ -v_{2} \end{smallmatrix} \right)}$. We compute the Euler class of $Y/\rho_{v}$. We denote the projection $Y \longrightarrow Y/\rho_{\left( \begin{smallmatrix} 1 \\ 0 \end{smallmatrix} \right)} \cong S_{\left( \begin{smallmatrix} 0 \\ 1 \end{smallmatrix} \right)}$ by $\pi_{\left( \begin{smallmatrix} 1 \\ 0 \end{smallmatrix} \right)}$ and the projection $Y \longrightarrow Y/\rho_{\left( \begin{smallmatrix} 0 \\ 1 \end{smallmatrix} \right)} \cong S_{\left( \begin{smallmatrix} 1 \\ 0 \end{smallmatrix} \right)}$ by $\pi_{\left( \begin{smallmatrix} 0 \\ 1 \end{smallmatrix} \right)}$. We fix an $S^1$-connection form $\eta_{v}$ on $S_{v}$ for $v=\left( \begin{smallmatrix} 1 \\ 0 \end{smallmatrix} \right)$ and $\left( \begin{smallmatrix} 0 \\ 1 \end{smallmatrix} \right)$. Let $\pi \colon Y \longrightarrow Y/\rho_{v}$ be the canonical projection. We take a $1$-form $\tilde{\eta}_{v}$ on $Y$ so that $\tilde{\eta}_{v}=v_{1} \pi_{\left( \begin{smallmatrix} 1 \\ 0 \end{smallmatrix} \right )}^{*}\eta_{\left( \begin{smallmatrix} 1 \\ 0 \end{smallmatrix} \right)} - v_{2} \pi_{\left( \begin{smallmatrix} 0 \\ 1 \end{smallmatrix} \right)}^{*} \eta_{\left( \begin{smallmatrix} 0 \\ 1 \end{smallmatrix} \right)}$ where we put $v=\left( \begin{smallmatrix} v_{1} \\ v_{2} \end{smallmatrix} \right)$. Then $\tilde{\eta}_{v}$ induces an $S^1$-connection form on $Y/\rho_{v}$. In fact, we have $\tilde{\eta}_{v}(\Xi(v))=0$ and $\tilde{\eta}_{v}(\Xi(w))=1$. Since $\tilde{\eta}_{v}$ is invariant by the $T^2$-action, $\tilde{\eta}_{v}$ is a basic form on the foliation defined by the fibers of $\pi$. Hence $\tilde{\tau}$ induces a $1$-form $\eta_{v}$ on $Y/\rho_{v}$. Because $\tilde{\eta}_{v}(\Xi(w))=1$, $\eta_{v}$ is an $S^1$-connection form. Since we have $\pi^{*}d\eta_{v}=v_{1} \pi_{\left( \begin{smallmatrix} 1 \\ 0 \end{smallmatrix} \right )}^{*} d\eta_{\left( \begin{smallmatrix} 1 \\ 0 \end{smallmatrix} \right)} - v_{2} \pi_{\left( \begin{smallmatrix} 0 \\ 1 \end{smallmatrix} \right)}^{*} d\eta_{\left( \begin{smallmatrix} 0 \\ 1 \end{smallmatrix} \right)}$, the Euler class of $Y/\rho_{v}$ is $v_{1}e(S_{\left( \begin{smallmatrix} 1 \\ 0 \end{smallmatrix} \right)}) - v_{2}e(S_{\left( \begin{smallmatrix} 0 \\ 1 \end{smallmatrix} \right)})=v_{1} [dx_{1} \wedge dx_{3}] - v_{2} [dx_{2} \wedge dx_{3}]$. Then $Y/\rho_{\left(\begin{smallmatrix} v_{1} \\ v_{2} \end{smallmatrix} \right)}$ is isomorphic to $S_{\left( \begin{smallmatrix} v_{1} \\ -v_{2} \end{smallmatrix} \right)}$.

Taking quotient of $Y$ by $\rho_{\left( \begin{smallmatrix} 2 \\ -1 \end{smallmatrix} \right)}$ and $\rho_{\left( \begin{smallmatrix} 1 \\ -1 \end{smallmatrix} \right)}$ respectively, $Y$ is written as
\begin{equation}\label{Equation : Y}
Y = \delta^{*} (\pr_{1}^{*} (Y/\rho_{\left( \begin{smallmatrix} 2 \\ -1 \end{smallmatrix} \right)}) \oplus \pr_{2}^{*} (Y/\rho_{\left( \begin{smallmatrix} 1 \\ -1 \end{smallmatrix} \right)})) = \delta^{*} ((\pr_{1}^{*} S_{\left( \begin{smallmatrix} 2 \\ 1 \end{smallmatrix} \right)}) \oplus (\pr_{2}^{*} S_{\left( \begin{smallmatrix} 1 \\ 1 \end{smallmatrix} \right)})).
\end{equation}

Let $f$ be the diffeomorphism on $T^4=\mathbb{R}^{4}/\mathbb{Z}^{4}$ induced by a matrix
\begin{equation}
A=\left( \begin{smallmatrix} 2 & 1 & 0 & 0 \\ 1 & 1 & 0 & 0 \\ 0 & 0 & 1 & 0 \\ 0 & 0 & 0 & 1 \end{smallmatrix} \right)
\end{equation}
in $GL(4;\mathbb{Z})$. Since $f$ commutes with $A$, $f$ induces a diffeomorphism $f_{1}$ on the orbifold $T^{4}/\sigma$. After blowing up singular points, $f_{1}$ extends to a diffeomorphism $f_{2}$ of the $K3$ surface $X$. The automorphism $f_{2}^{*}$ of $H^{2}(X;\mathbb{R})$ induced by $f_{2}$ has an invariant space $V=\mathbb{R} [dx_{1} \wedge dx_{3}] \oplus \mathbb{R} [dx_{2} \wedge dx_{3}]$. By the definition of $f_{2}$, we have
\begin{equation}\label{Equation : S}
f_{2}^{*}S_{\left( \begin{smallmatrix} 2 \\ 1 \end{smallmatrix} \right)} \cong S_{\left( \begin{smallmatrix} 1 \\ 0 \end{smallmatrix} \right)}, \quad f_{2}^{*}S_{\left( \begin{smallmatrix} 1 \\ 1 \end{smallmatrix} \right)} \cong S_{\left( \begin{smallmatrix} 0 \\ 1 \end{smallmatrix} \right)}.
\end{equation}
Hence we have $f_{2}^{*}Y = Y$ by \eqref{Equation : Y}. Then $f_{2}$ induces a diffeomorphism $\tilde{f}_{2}$ on $Y$. By \eqref{Equation : S}, we have
\begin{equation}\label{Equation : S 2}
\tilde{f}_{2 *}(\Xi{\left( \begin{smallmatrix} 1 \\ 0 \end{smallmatrix} \right)}) =\Xi \left( \begin{smallmatrix} 2 \\ 1 \end{smallmatrix} \right), \quad \tilde{f}_{2 *}(\Xi{\left( \begin{smallmatrix} 0 \\ 1 \end{smallmatrix} \right)}) = \Xi \left( \begin{smallmatrix} 1 \\ 1 \end{smallmatrix} \right).
\end{equation}
Then $\tilde{f}_{2 *}\Xi \left( \begin{smallmatrix} 1 \\ \frac{- 1 - \sqrt{5}}{2} \end{smallmatrix}\right) =\Xi \left( \begin{smallmatrix} 1 \\ \frac{- 1 - \sqrt{5}}{2} \end{smallmatrix}\right)$. Hence $\tilde{f}_{2}$ preserves $\mathcal{G}$.

Let $M$ be the mapping torus $([0,1] \times Y)/(0,\tilde{f}_{2}(x)) \sim (1,x)$. Let $\mathcal{F}$ be the foliation on $M$ induced from the foliation $\mathcal{F}_{0}$ with leaves $\{t\} \times L$, for $t \in [0,1]$ and $L \in \mathcal{G}$, on $[0,1] \times Y$. $\mathcal{F}$ is a Riemannian foliation, because $(M,\mathcal{F})$ has a bundle-like metric induced by the bundle-like metric $r(t)g + (1-r(t)) \tilde{f}_{2}^{*}g + dt \otimes dt$ on $([0,1] \times Y, \mathcal{F}_{0})$, where $g$ is a bundle-like metric of $(Y,\mathcal{G})$ and $r(t)$ is a smooth non-negative function on $[0,1]$ such that $r(t)=0$ near $0$ and $r(t)=1$ near $1$.

We show that $(M,\mathcal{F})$ is not minimizable. Let $\xi(\mathcal{F})$ be the \'{A}lvarez class of $(M,\mathcal{F})$. Since $(N,\mathcal{G})$ is an isometric flow, $(N,\mathcal{G})$ is minimizable. Hence the restriction of $\xi(\mathcal{F})$ to any fiber of the projection $p \colon M \longrightarrow S^{1}$ defined by $p([t,x])=[t]$ is zero. Let $\gamma$ be a closed path in $M$ which is mapped to a generator of $\pi_{1}S^1$. Since the closures of leaves of $(M,\mathcal{F})$ are $T^2$ with a linear flow of slope $\frac{- 1 - \sqrt{5}}{2}$, the structural Lie algebra of $(M,\mathcal{F})$ is $\mathbb{R}$. Since $\tilde{f}_{2}$ acts to the normal bundle of the linear flows on the closures of leaves of $\mathcal{F}$ by the multiplication of $\frac{3 - \sqrt{5}}{2}$, the holonomy of the Molino's commuting sheaf of $(M,\mathcal{F})$ with respect to $\gamma$ is the map defined by the multiplication of $\frac{3 - \sqrt{5}}{2}$ by the diagram \eqref{Equation : Diagram}. By a theorem of \'{A}lvarez L\'{o}pez \cite{Alvarez Lopez 2}, we have
\begin{equation}
\int_{\gamma}\xi(\mathcal{F}) = \log \frac{3 - \sqrt{5}}{2}.
\end{equation}
Hence, by a theorem of \'{A}lvarez L\'{o}pez \cite{Alvarez Lopez}, $(M,\mathcal{F})$ is not minimizable.

By the homotopy exact sequence of a principal $S^1$-bundle $E$ over a simply connected manifold $B$, $\pi_{1}E$ is isomorphic to a quotient of $\mathbb{Z}$. Then $\pi_{1}E$ is generated by a closed path along an $S^1$-fiber of $E$. Hence $E$ is simply connected if and only if there exists no covering space of $E$ whose restriction to the fiber of an $S^1$-fiber is finite and nontrivial. If $E'$ is a principal $S^1$-bundle and an $n$-fold covering space of $E$ along the $S^1$-fibers, we have $n e(E') = e(E)$. Hence $e(E)$ is divisible by $n$. Thus $Y$ is simply connected by construction. Hence $\pi_{1}M$ is isomorphic to $\mathbb{Z}$ by the homotopy exact sequence of the fiber bundle $Y \longrightarrow M \longrightarrow S^1$.

\subsection{Riemannian foliations of codimension $2$ on fiber bundles over $S^{1}$}

Let $(N,\mathcal{G})$ be a Riemannian foliation of codimension $1$ on a closed manifold $N$. Let $f$ be a diffeomorphism on $N$ which preserves $\mathcal{G}$. Let $M=([0,1] \times N)/(0,f(x)) \sim (1,x)$ be the mapping torus of $f$. Let $\mathcal{F}$ be a foliation on $M$ induced from the foliation on $[0,1] \times N$ with leaves $\{t\} \times L$, for $t \in [0,1]$ and $L \in \mathcal{G}$. $\mathcal{F}$ is Riemannian by the same reason as the previous example. If $\pi_{1}M$ is of polynomial growth, then $(M,\mathcal{F})$ is minimizable by Corollary \ref{Corollary : Codimension 2}.

\end{document}